\documentclass[preprint,12pt]{elsarticle}

\newif\iffastcompile
\ifdefined\FASTCOMPILE
	\fastcompiletrue
\else
	\fastcompilefalse
\fi
\iffastcompile
	\PassOptionsToPackage{draft}{graphicx}
\fi

\usepackage{amssymb}
\usepackage{amsmath}

\usepackage{siunitx}
\usepackage{booktabs}
\usepackage{url}
\iffastcompile
	\usepackage[draft,bookmarks=false]{hyperref}
\else
	\usepackage[bookmarks=true]{hyperref}
\fi
\iffastcompile
	\setkeys{Gin}{draft=true}
\fi

\journal{Advannces in Space Research}

\begin{document}

\begin{frontmatter}



  \title{Fuel-Optimal Low-Thrust Trajectory Design under High-Fidelity Dynamics: A State Transition Matrix-Based Sensitivity Approach}

  \author{Liqiang Hou}
  \address{houliqiang2008@gmail.com, Shanghai, China}

  \begin{abstract}
    A straightforward and computationally efficient indirect method based on STM sensitivity analysis is introduced for the design of fuel-optimal low-thrust transfers under high-fidelity dynamics. In conventional indirect approaches, explicit expressions for the partial derivatives of the system dynamics are required to formulate the costate equations. Deriving these explicit costate equations is complex and challenging for high-fidelity trajectory design. In this study, the costate equations of the optimal control problem are reformulated as a set of ordinary differential equations involving the state variables and their time derivatives. Complex dynamical effects are treated as black-box components in the design, thereby avoiding tedious derivations of the partial derivatives of the system dynamics. A standard gradient or interior-point based optimizer is employed to determine the optimal costates and transfer parameters. The equivalence between the proposed method and conventional approach is demonstrated through a classic Earth--Mars transfer simulation scenario. An Earth--Mars transfer design under high-fidelity dynamics is then presented using the proposed method. Perturbations due to solar radiation pressure, solar $J_2$ oblateness, Jupiter’s third-body gravity, and relativistic effects are considered. Finally, application of the method to a multiple revolution Earth-Venus transfer under high-fidelity dynamics is  presented and discussed.
  \end{abstract}

  \begin{graphicalabstract}
  \end{graphicalabstract}


  \begin{highlights}
    \item A new indirect approach for low-thrust transfer design under high-fidelity dynamics.
    \item New reformulated costate-equation expression for control-law design, avoiding tedious derivations of system-dynamics partial derivatives.
    \item Complex dynamical effects are treated as black-box components in the design process.
    \item Equivalence between the proposed method and conventional approaches is demonstrated through classic transfer simulation scenario.
    \item Performance of the method is demonstrated through different mission scenarios and transfer targets.
  \end{highlights}

  \begin{keyword}
    Design optimization \sep State Transition Matrix \sep Low-thrust transfer \sep High-fidelity dynamics \sep Costate analysis
  \end{keyword}

\end{frontmatter}



\section{Introduction}
\label{sec: Introduction}

Low-thrust, fuel-optimal transfer design in space mission analysis is typically formulated as an initial-value problem governed by a set of ordinary differential equations (ODEs), including the equations of motion and the associated costate equations \cite{Morante2021}. The costate equations are expressed through the partial derivatives of the system dynamics. In general, these equations are difficult to derive and write explicitly in closed form. In high-fidelity mission design, additional perturbative effects such as higher-order gravitational harmonics, solar radiation pressure, and relativistic corrections must be considered. These effects introduce strong nonlinearities, and their partial derivatives with respect to the state variables are often analytically intractable.

In cases with comparatively simple dynamics, such as two-body and three-body gravity models, explicit expressions for the costate equations can be derived \cite{zuiani2020multiobjective, leomanni2019indirect}. The optimal control law for fuel-optimal transfer can then be obtained by solving a two-point boundary value problem (TPBVP), in which the equations of motion and the explicit costate equations are integrated to satisfy the terminal equality constraints \cite{jiang2016adjoint}. The resulting control law typically consists of a sequence of thrust-on and thrust-off phases. Because of these switching characteristics, standard gradient-based optimization techniques are often insufficient for accurately determining the initial costates. Indirect approaches, such as continuation or homotopy methods, are therefore frequently employed, in which a simpler problem, such as an energy-optimal transfer, is used as a starting point and gradually transformed into the fuel-optimal solution \cite{taheri2017homotopic,taheri2018continuation}.

Direct optimization methods, including collocation and shape-based parameterization approaches \cite{Morante2021, izzo2021trajectory,PattersonRao2014GPOPS}, can also be used for transfer design. These methods discretize the state and control variables, allowing high-fidelity dynamical models to be incorporated without requiring explicit derivation of the costate equations. The accuracy and convergence behavior of these methods are influenced by the chosen collocation strategy and the distribution of discretization points, particularly in large-scale or highly complex design scenarios. Improperly designed discretization schemes can lead to significant numerical errors. Other methods based on hybrid approaches \cite{englander2017automated} have also been employed for trajectory optimization problems. Recent studies have also proposed learning-based techniques for low-thrust trajectory optimization.

In this study, a new framework for the costate determination of fuel-optimal transfer design is proposed. The state transition matrix (STM), together with sensitivity-based analysis, is incorporated into the optimization framework. The STM technique had been employed in impulsive transfer scenarios, including relative-motion reconfiguration and transfer design under restricted three-body dynamics \cite{koenig2017stm,howell2020halo,koenig2016relative}. In this study, a new STM-based sensitivity analysis is developed for low-thrust costate equations. New differential equations for the costate are derived and expressed in terms of the state variables and their time derivatives, rather than the conventional partial derivatives of the system dynamics, thereby avoiding tedious derivations of partial differential expressions. High-fidelity perturbations are incorporated directly into these equations and subsequently integrated into the optimization of the control law, which is essential for design under complex dynamical perturbations.

A classical fuel-optimal design scenario, the low-thrust Earth-Mars transfer, is used to validate the proposed method. First, a trajectory design under two-body orbital dynamics is simulated. The initial costate solutions are computed and compared with those obtained from the classical costate-based formulation using explicit partial derivatives, demonstrating equivalence of the two approaches. Subsequently, a high-fidelity low-thrust transfer using the proposed method is presented. Perturbation impacts due to J2 gravitational harmonics, solar radiation pressure, third-body gravity from Jupiter, and relativistic corrections on the Earth-Mars transfer, etc., are considered. Simulation results for the high-fidelity model are listed and compared to those from the low-fidelity model. To further illustrate the method, a multiple-revolution Earth-Venus transfer under low-thrust control, involving multiple coast and thrust arcs, is also simulated.

The remainder of the paper is organized as follows. Section \ref{sec:model E-M} briefly reviews the optimal-control formulation for low-thrust transfer design. Equations of motion and explicit expressions for costate differential equations under two-body dynamics, are listed. The proposed method, based on state transition matrix  and sensitivity analysis, is presented in Section \ref{sec:State transition and sensitivity analysis}. Transfer design under high-fidelity dynamics is then developed and presented.  Numerical simulations and analysis of the method, including Earth-Mars transfers under low- and high-fidelity dynamics, are presented in Section \ref{sec:Numerical Simulation}. Equivalence test between the proposed method and the conventional approach is conducted, along with an application to a multiple revolution Earth-Venus transfer. Finally, conclusions and future work are summarized in Section \ref{sec:conclusion}.

\section{Fuel-Optimal Low-Thrust Control-Law Design}
\label{sec:model E-M}

\subsection{Equations of Motion: the Low-Fidelity Model}

The dynamics of the space transfer is modeled and given as (for Cartesian mode)

\begin{equation}
  \begin{aligned}
    \dot{\mathbf{r}} & =\mathbf{v},                                                 \\
    \dot{\mathbf{v}} & =\mathbf{g}(\mathbf{r})+c_1 \frac{u}{m} \boldsymbol{\alpha}, \\
    \dot{m}          & =-c_2 u \label{eqn:equation of motion}
  \end{aligned}
\end{equation}
where $\mathbf{r}$ and $\mathbf{v}$ are  position and velocity of the spacecraft. $u \in [0, 1]$, $\boldsymbol{\alpha} \in [0,1]^{3}$ are the thrust ratio and vector to be optimized. Constants of the thrust,  $c_1$ and $c_2$, are defined as
\begin{equation}
  c_1  =T_{\max }, c_2 =\frac{T_{\max }}{I_{s p} g_0}
\end{equation}
where $T_{\text {max }}$ is the maximal thrust magnitude, $I_{s p}$ is the thruster specific impulse and $g_0$ is the standard acceleration of gravity at sea level, 9.80665 m/$s^2$. $\mathbf{g}(\mathbf{r})$ is gravity function. For transfer design under two-body dynamics, the gravity function is set to
\begin{equation}\label{eqn:two-body}
  \mathbf{g}(\mathbf{r})=-\frac{\mu_{s}}{r^3} \mathbf{r},
\end{equation}
where $\mu_{s}$ is the Sun gravitational constant, and set to the constant 1.0 in canonical heliocentric units.

\subsection{Equations of Motion: the High-Fidelity Model}

Consider a space trajectory design under high order  heliocentric acceleration computations. In addition to central Sun gravity, perturbations from Jupiter (dominant third body), solar radiation pressure, solar J2 and first-order post-Newtonian relativistic correction, etc., are considered.  The total acceleration of dynamics consists of
\begin{equation}\label{eqn:high_fidelity_acc}
  \mathbf{a}=\mathbf{a}_s+\sum_i \mathbf{a}_i+\mathbf{a}_{J 2}+\mathbf{a}_{S R P}+\mathbf{a}_{G R}
\end{equation}
where $\mathbf{a}_s$, $\mathbf{a}_i$, $\mathbf{a}_{J 2}$, $\mathbf{a}_{S R P}$, and $\mathbf{a}_{G R}$ are accelaertaions due to the central solar gravity, third-body perturbations, solar oblateness (J2), solar radiation pressure (SRP), and  relativistic correction, respectively. Introduced high order acceleration items are briefly summarized  as below.

The dominant acceleration for the heliocentric transfer design, $\mathbf{a}_s$, is computed as
\begin{equation*}
  \mathbf{a}_{s}=-\frac{\mu_{s}}{r^3} \mathbf{r}
\end{equation*}
Acceleration due to the solar oblateness, $\mathbf{a}_{J 2}$, is computed as
\begin{equation*}
  \mathbf{a}_{J 2}=\frac{3 J_2 \mu_{s} R_{\odot}^2}{2 r^5}\left[\begin{array}{l}
      x\left(5 z^2 / r^2-1\right) \\
      y\left(5 z^2 / r^2-1\right) \\
      z\left(5 z^2 / r^2-3\right)
    \end{array}\right]
\end{equation*}
where $J_2$ is solar oblateness coefficient, and $R_{\odot}$ is solar radius. Perturbation due to planet $i$ is computed as
\begin{equation*}
  \mathbf{a}_i=\mu_i\left(\frac{\mathbf{r}_i-\mathbf{r}}{\left|\mathbf{r}_i-\mathbf{r}\right|^3}-\frac{\mathbf{r}_i}{\left|\mathbf{r}_i\right|^3}\right)
\end{equation*}
where $\mathbf{r}_i$ is the heliocentric planet position, and $\mu_i=G M_i$ is the planet's gravity constant.  The planets in this study are set to $i = \{\text{Venus}, \text{Earth}, \text{Mars}, \text{Jupiter}\}$. Among the perturbations, long-period perturbation due to the Jupiter is particularly important, since $\mu_J \approx 9.545 \times 10^{-4} \mu_{s}$.

The solar radiation in the transfer varies closely with the area-to-mass ratio, leads to small yet long-duration orbit and attitude deviations. Acceleration due to the solar radiation pressure is computed as
\begin{equation*}
  \mathbf{a}_{S R P}=\beta \frac{\mu_{s}}{r^3} \mathbf{r}
\end{equation*}
where $\beta=\frac{C_r A P_0}{m}$.  $C_r$, $\frac{A}{m}$ and $P_0$ are reflectivity coefficient, area-to-mass ratio and solar radiation pressure at 1 AU, respectively. The Schwarzschild correction for the general relativity is
\begin{equation*}
  \mathbf{a}_{G R}=\frac{\mu_{s}}{c^2 r^3}\left[\left(\frac{4 \mu_{s}}{r}-v^2\right) \mathbf{r}+4(\mathbf{r} \cdot \mathbf{v}) \mathbf{v}\right]
\end{equation*}
where $c=$ speed of light and $v=|\mathbf{v}|$.

Take into account total acceleration due to above factors, 	equation of the motion under high fidelity model can be listed as

\begin{equation}\label{eqn:high_fidelity_propagator}
  \begin{aligned}
    \dot{\mathbf{r}} & =\mathbf{v},                                                                                                                                                                                                                                                                 \\
    \dot{\mathbf{v}} & = -\frac{\mu_{s}}{r^3} \mathbf{r}+\sum_i \mu_i\left(\frac{\mathbf{r}_i-\mathbf{r}}{\left|\mathbf{r}_i-\mathbf{r}\right|^3}-\frac{\mathbf{r}_i}{\left|\mathbf{r}_i\right|^3}\right)+\mathbf{a}_{J 2}+\mathbf{a}_{S R P}+\mathbf{a}_{G R} +c_1 \frac{u}{m} \boldsymbol{\alpha} \\
    \dot{m}          & =-c_2 u
  \end{aligned}
\end{equation}
where $c_1$, $u$ and $\boldsymbol{\alpha}$ are thrust related parameters as mentioned in preceding sections.

\subsection{The Optimal Control Formulation for Fuel-Optimal Space Trajectory Design}

Consider a fuel-optimal space trajectory design. Differential equations for the low-fidelity and high-fidelity models,  can be written as
\begin{equation}
  \dot{\mathbf{x}} = \mathbf{f}(\mathbf{x} , \mathbf{u}, t)
\end{equation}
with state variables $\mathbf{x} = [\mathbf{r},\mathbf{v},m]$, and control variable $\mathbf{u}$. The control variable $\mathbf{u}$ is set to the thrust direction $\boldsymbol{\alpha}$ and the thrust ratio $u$
\begin{equation}
  \mathbf{u}=u \boldsymbol{\alpha} .
\end{equation}
with the constraints
\begin{equation}
  \begin{array}{r}
    u \leq 1, \\
    \boldsymbol{\alpha} \cdot \boldsymbol{\alpha}-1=0 .
  \end{array}
\end{equation}
The time-dependent objective of fuel consumption is
\begin{equation}\label{eqn: objective}
  \mathcal{J}=\int_{t_0}^{t_f} \mathcal{L}(\mathbf{x}, \mathbf{u}, t) d t =  \int_{t_0}^{t_f} c_2 u(t) d t
\end{equation}
with the constraints at final time $t_f$
\begin{equation}
  \mathbf{\psi}^f\left(\mathbf{x}\left(t_f\right), t_f\right)=0
\end{equation}
For a rendezvous mission design, the terminal constraint is defined as $\mathbf{\psi}^f:=[\mathbf{r}\left(t_f\right)-\bar{\mathbf{r}}_f, \mathbf{v}\left(t_f\right)-\bar{\mathbf{v}}_f]$, where $\bar{\mathbf{r}}_f$ and $\bar{\mathbf{v}}_f$ denote the desired terminal position and velocity of the target planet, respectively.

Adjoining the system differential equations and the terminal constraints, new equivalent unconstrained objective can be defined using the Lagrange multiplier functions and costate of the system equation, $\boldsymbol{\lambda}(t)$
\begin{equation}\label{eq:performance index}
  \mathcal{J}= \mathbf{\psi}^f\left(\mathbf{x}\left(t_f\right), t_f\right)+\int_{t_0}^{t_f}\left\{\mathcal{L}(\mathbf{x}, \mathbf{u}, t)+\boldsymbol{\lambda}^T(\mathbf{f}(\mathbf{x}, \mathbf{u}, t)-\dot{\mathbf{x}})\right\} d t
\end{equation}
Define Hamiltonian of the system
\begin{equation}
  \mathcal{H}=\mathcal{L}+\boldsymbol{\lambda}^T \mathbf{f}=\boldsymbol{\lambda}_{\boldsymbol{r}} \cdot \boldsymbol{v}+\boldsymbol{\lambda}_{\boldsymbol{v}} \cdot\left(\mathbf{g}(\mathbf{r})+c_1 \frac{u}{m} \boldsymbol{\alpha}\right)-\lambda_m c_2 u+c_2 u
\end{equation}
the equation becomes
\begin{equation}
  \mathcal{J}=\mathbf{\nu}^T \mathbf{\psi}^f\left(\mathbf{x}\left(t_f\right), t_f\right)+\int_{t_0}^{t_f}\left\{\mathcal{H}-\boldsymbol{\lambda}^T \dot{\mathbf{x}}\right\} d t \text {. }
\end{equation}

The necessary optimality conditions for the control are derived from Pontryagin’s  principle. \cite{pontryagin2018mathematical}. According to the principle, the admissible controls are selected such that the Hamiltonian $\mathcal{H}$ is minimized at each point along the trajectory. Consequently, the optimal thrust direction is determined by the velocity costate $\boldsymbol{\lambda}_v$
\begin{equation}\label{eqn: optimal alpha}
  \boldsymbol{\alpha}=-\boldsymbol{\lambda}_v / \lambda_v,
\end{equation}
with the primer vector-based control law

\begin{equation}\label{eqn: primer vector control law}
  u= \begin{cases}0 & \text { if } \quad S>0 \\ 1 & \text { if } \quad S<0 \\ 0 \leq u \leq 1 & \text { if } \quad S=0 .\end{cases}
\end{equation}
where $S$ is the switching function
\begin{equation}\label{eqn: switch function}
  S=1-\frac{c_1}{c_2 m} \lambda_v-\lambda_m
\end{equation}

\subsection{Costate Equations and the TPBVP}
Differential equation of the costate $\boldsymbol\lambda$ is defined as $\dot{\boldsymbol{\lambda}}  =-\mathcal{H}_x$. For fuel-optimal transfer design, the costate equations can be expressed in terms of partial derivatives of the system dynamics
\begin{equation}\label{eqn:costate with Jacobian}
  \dot{\boldsymbol{\lambda}}=-\frac{\partial \mathcal{H}}{\partial \mathbf{x}}=- \left(\frac{\partial \mathbf{f}}{\partial \mathbf{x}}\right)^T \boldsymbol{\lambda}
\end{equation}
Transversality conditions of the design are defined as  $\boldsymbol{\lambda}^T\left(t_f\right)  = \mathbf{\psi}_x^f$. The final mass costate is zero according to the transversality condition, since the final mass is free
\begin{equation}\label{eqn:lambda mass}
  \lambda_m\left(t_f\right)=0
\end{equation}

Incorporating the transversality conditions and the terminal constraints, TPBVP for the optimal control problem can be formulated as
\begin{equation*}
  \boldsymbol\Psi^0: \mathbf{r}\left(t_0\right)=\mathbf{r}_0, \quad \mathbf{v}\left(t_0\right)=\mathbf{v}_0, \quad m\left(t_0\right)=1,
\end{equation*}
and
\begin{equation*}
  \boldsymbol\Psi^f:	\mathbf{r}\left(t_f\right)=\bar{\mathbf{r}}_f, \quad \mathbf{v}\left(t_f\right)=\bar{\mathbf{v}}_f, \quad \lambda_{m,f} = 0
\end{equation*}
The design variables for solving the TPBVP by shooting are the initial costate values, $\boldsymbol{\lambda}_0 = [\boldsymbol{\lambda}_{r,0},\boldsymbol{\lambda}_{v,0},\lambda_{m,0}]$.

In the costate equation, Eq.~\eqref{eqn:costate with Jacobian}, Jacobian of the system dynamics is required. For comparatively simple cases, such as trajectory design under two-body or three-body gravity models, explicit expressions for the costate ODEs can be derived. Equation~\eqref{eqn: two body costate} presents explicit costate equations for the two-body gravitational model.
\begin{equation}\label{eqn: two body costate}
  \begin{aligned}
    \dot{\boldsymbol{\lambda}}_r & =\frac{\mu}{r^3} \boldsymbol{\lambda}_v-\frac{3 \mu \mathbf{r} \cdot \boldsymbol{\lambda}_v}{r^5} \mathbf{r} \\
    \dot{\boldsymbol{\lambda}}_v & =-\boldsymbol{\lambda}_r                                                                                     \\
    \dot{\lambda}_m              & =-c_1 \frac{u}{m^2} \lambda_v .
  \end{aligned}
\end{equation}
The TPBVP is highly sensitive to the initial guesses of the costates due to the nonlinear effects of the switching function. Shooting methods based on gradient-search algorithms, such as Newton’s or Powell’s method, often struggle to solve this problem, even for two-body gravitational models. Indirect strategies, such as continuation methods or homotopy smoothing, are therefore commonly employed. As for complex scenarios involving higher-order nonlinear perturbations, obtaining an explicit analytical expression for the costate equations is mathematically challenging, even with advanced symbolic computation tools.

\section{State Transition and Sensitivity Analysis for Transfer Design}
\label{sec:State transition and sensitivity analysis}

To address the challenges arising from the PDE-like nature of the associated costate equations, a numerical STM-based costate analysis is developed for control design. Variations in the initial state are propagated through the variational equations governed by the system Jacobian. Based on this formulation, a state-transition-matrix approach is then established for costate analysis.

\subsection{State Transition Matrix (STM) }
Given the dynamics system with the flow defined as $\mathbf{x}=\boldsymbol{\phi}\left(t; \mathbf{x}_0, t_0\right)$, a small perturbation on a reference solution of the flow can be written as
\begin{equation*}
  \mathbf{x}+\delta \mathbf{x}=\boldsymbol{\phi}\left(t;\mathbf{x}_0+\delta \mathbf{x}_0, t_0 \right)
\end{equation*}
Corresponding first-order expansion  is given as
\begin{equation*}
  \mathbf{x}+\delta \mathbf{x}=\boldsymbol{\phi}\left(t;\mathbf{x}_0, t_0\right)+\frac{\partial \boldsymbol{\phi}}{\partial \mathbf{x}_0}\left(t;\mathbf{x}_0, t_0\right) \delta \mathbf{x}_0
\end{equation*}
With the equations, relation of a perturbation at $t$ with respect to the perturbation of the initial state is found as,

\begin{equation*}
  \delta \mathbf{x}=\boldsymbol{\Phi}\left(t, t_0\right) \delta \mathbf{x}_0
\end{equation*}
where $\boldsymbol{\Phi}\left(t, t_0\right)$ indicates the transition matrix defined as,

\begin{equation*}
  \boldsymbol{\Phi}\left(t, t_0\right)=\frac{\partial \boldsymbol{\phi}}{\partial \mathbf{x}_0}\left(t;\mathbf{x}_0, t_0\right)
\end{equation*}

Extend the analysis to the trajectory design of nonlinear dynamics $ \dot{\mathbf{x}}=\mathbf{f}(\mathbf{x}, t)$. STM and its ODE of the dynamics can be given as
\begin{equation}\label{eqn:STM dynamics}
  \dot{\mathbf{\Phi}}\left(t, t_0\right)=\mathbf{A}(t) \mathbf{\Phi}\left(t, t_0\right)
\end{equation}
where $\mathbf{A}(t)$ is the Jacobian of the system. In the two-body dynamics, the matrix $\mathbf{A}(t)$ can be expressed as \cite{koenig2017stm},
\begin{equation}
  \mathbf{A}(t)=\left[\begin{array}{cc}
                 0_{3 \times 3}     & \mathbf{I}_{3 \times 3} \\
                 \mathbf{U}_{X X}^* & 2 \mathbf{\Omega}
    \end{array}\right],
\end{equation}
where $\mathbf{U}_{X X}^*$ is the Hessian of the potential function at the reference  $\mathbf{x}^*$, and $\mathbf{\Omega}$ is the cross product matrix of the angular velocity of the rotating frame.

The STM maps the state at time $t$ relative to the state at the initial epoch $t_0$. At $t_0$, the STM is the identity matrix, i.e., $\boldsymbol{\Phi}\left(t_0, t_0\right) = \mathbf{I}$. Having the initial conditions well defined, the STM can be propagated by utilizing the system dynamics. Usages of STM in space designs include multiple-impulse transfers, stability analysis, and manifold-based computations, etc \cite{koenig2017stm,englander2017automated}.

\subsection{STM-Based Reformulation of the Costate Equation}
Consider a determination of the costates for the low-thrust space trajectory control law. The trajectory state can be implicitly propagated by integrating the coupled state and costate ODEs. Given an initial guess of $\boldsymbol{\lambda}_{0}$, this propagation process can be expressed as the nonlinear dynamical system $\dot{\mathbf{x}}=\boldsymbol{f}(\mathbf{x},t;\boldsymbol{\lambda}_0)$. The optimization problem is then to determine an appropriate $\boldsymbol{\lambda}_{0}$ such that the terminal-state and terminal-costate constraints are satisfied within a predefined tolerance.

Expanding the trajectory to first order about the reference state $\mathbf{x}_0$, and following the STM analysis in the preceding section, dynamics of the state displacement, $\delta \mathbf{x}(t)$, can be written as
\begin{equation}
  \begin{aligned}
    \delta \mathbf{x}              & =\boldsymbol{\Phi}(t,t_0) \delta \mathbf{x}_{0} \\
    \dot{\boldsymbol{\Phi}}(t,t_0) & = \mathbf{A}(t) \Phi\left(t, t_{0}\right)
  \end{aligned}
\end{equation}
where $\mathbf{A}$ is Jacobian of the system
\begin{equation}
  \mathbf{A}=\frac{\partial \mathbf{f(\mathbf{x},t;\boldsymbol{\lambda}_0)}}{\partial \mathbf{x}}
\end{equation}
In this analysis, the state at time $t_0$, $\mathbf{x}(t_0;\boldsymbol{\lambda}_0)$, is not necessarily equal to the prescribed initial state $\mathbf{x}_0$, due to the introduced costate $\boldsymbol{\lambda}_0$ and associated control $\boldsymbol{u}_0$.

The costate equation for the fuel-optimal transfer in Eq.\eqref{eqn:costate with Jacobian}, assuming constant thrust and specific impulse, can be written as
$\dot{\boldsymbol{\lambda}} =  -\mathbf{A}^T \boldsymbol{\lambda}$.
Multiplying both sides of the equation by the outer product of the state deviation, $\delta \mathbf{x}\delta \mathbf{x}^T$, yields
\begin{equation}\label{eqn:stm-costat-step-1}
  \begin{aligned}
    (\delta \mathbf{x}\cdot \delta \mathbf{x}^T)  \dot{\boldsymbol{\lambda}}
     & = -(\delta \mathbf{x}\cdot \delta \mathbf{x}^T)  \mathbf{A}^T \boldsymbol{\lambda}                                         \\
     & =  -\delta \mathbf{x}\cdot \left(  \boldsymbol{\Phi}(t,t_0)\delta \mathbf{x}_0 \right) ^T \mathbf{A}^T\boldsymbol{\lambda} \\
     & =  -\delta \mathbf{x}\cdot \delta \dot{\mathbf{x}}^T \boldsymbol{\lambda}                                                  \\
  \end{aligned}
\end{equation}
In the equation, prediction of the time derivative of the deviation, $\delta\dot{ \mathbf{x}}$ , is obatined through STM propagation
\begin{equation}
  \delta\dot{ \mathbf{x}} =	\mathbf{A}  \boldsymbol{\Phi}(t,t_0)\delta \mathbf{x}_0
\end{equation}

In Eq.~\eqref{eqn:stm-costat-step-1}, the outer products $\delta \mathbf{x}\delta \mathbf{x}^T$ and $\delta \mathbf{x}\delta \dot{\mathbf{x}}^T$ are both rank-one matrices; therefore, their inverses are singular. By multiplying Eq.~\eqref{eqn:stm-costat-step-1} again by $\delta \mathbf{x}\delta \mathbf{x}^T$ and rearranging the resulting terms, one can have the STM-based costate equation
\begin{equation}\label{eqn:STM-based costate equation}
  \dot{\boldsymbol{\lambda}}  = -\left(\delta \mathbf{x}\delta \mathbf{x}^T + \epsilon \mathbf{I} \right)^{+} \delta \mathbf{x} \delta \dot{\mathbf{x}}^T   \boldsymbol{\lambda} \end{equation}
where $\left(\delta \mathbf{x}\delta \mathbf{x}^T + \epsilon \mathbf{I}\right)^{+}$ denotes the pseudoinverse of the regularized outer product $\delta \mathbf{x}\delta \mathbf{x}^T$, and $\epsilon$ is a sufficiently small positive constant introduced to avoid singularity in the computation.

\begin{figure}[h]
  \centering
  \includegraphics[width=0.7\linewidth]{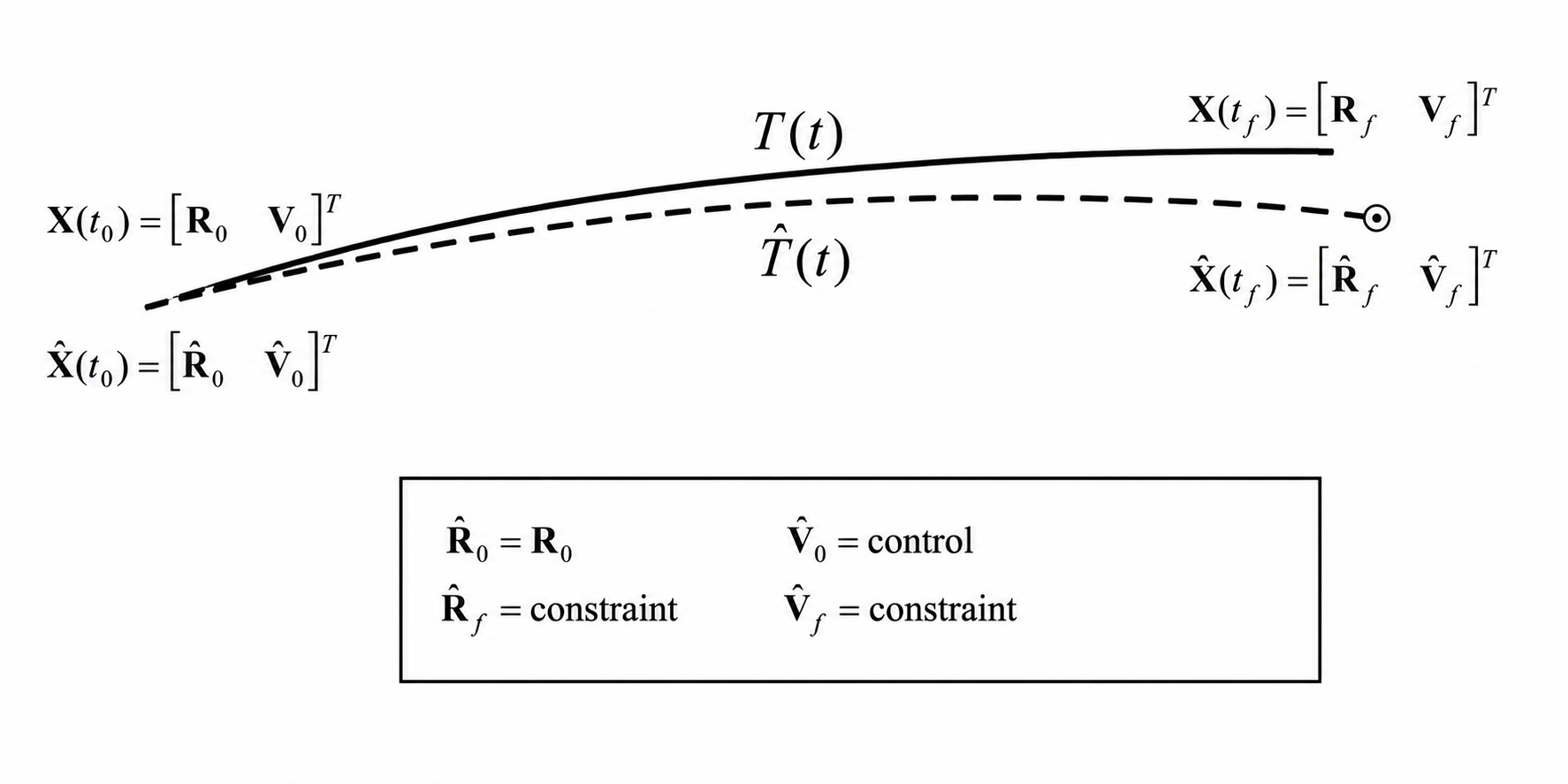}
  \caption{Differential-correction shooting routine for the space transfer. The solid and dashed lines denote the initial trajectory and the corrected trajectory that reaches the target position at the target time, respectively. The forward trajectory is propagated using the current initial guess and costate equation; the initial costates are then corrected using STM-based sensitivity information, and the process is repeated until the terminal constraints are satisfied.}
  \label{fig:the-single-shooting-differential-correction-routine}
\end{figure}

The new costate equation, together with the state equation and the terminal constraints $\boldsymbol{\Psi}_f$, forms a TPBVP for the transfer design. The switching function, Eq.~\eqref{eqn: switch function}, and the primer vector in the control law, Eq.~\eqref{eqn: primer vector control law}, remain unchanged because the costate equation is derived from the same Pontryagin principles as the original formulation. Figure~\ref{fig:the-single-shooting-differential-correction-routine} schematically illustrates the shooting routine for the STM-based space transfer design.

In the new STM-based costate computation, the system Jacobian is replaced by numerical computations using the state variables and their time derivatives. The resulting equations can therefore be directly applied to transfer designs under high-order system models, without requiring an explicit mathematical expression for the system Jacobian in the costate computation. This feature is particularly useful for those designs under complex dynamics and high-fidelity models. Computational performances of the proposed method are demonstrated through a set of space trajectory design simulations.

\subsection{Equivalence to the Original Costate Equation in Two-Body Transfer Design}

It is known that the solution of the fuel-optimal costate is not unique. To demonstrate equivalence of the proposed costate computation to the original formulation in Eq.~\eqref{eqn:costate with Jacobian}, a two-body transfer design for which the costate equation can be expressed explicitly is considered. An augmented design parameter, $\nu_2 > 0$, is introduced into the fuel objective function
\begin{equation}\label{eqn:obj equation with nu}
  J_{\nu} = \int_{t_0}^{t_f} \mathcal{L}_{\nu} dt= \int_{t_0}^{t_f}\left(\nu_2 c_2 u \right) dt
\end{equation}
The parameter $\nu_2$ provides additional freedom in determining the initial costates while preserving the optimality of the design.

The terminal constraint on the position and velocity of the transfer, given the reference $\mathbf{r}_0$ and $\mathbf{v}_0$ at the departure, can be written equivalently as
\begin{equation}
  \begin{aligned}
    \boldsymbol{\psi}_{\nu}^f : \quad & \delta \mathbf{r}_\nu^f = (\mathbf{r}(t_f) - \mathbf{r}_0) - (\bar{\mathbf{r}}_f - x_0), \\
                                      & \delta \mathbf{v}_\nu^f = (\mathbf{v}(t_f) - \mathbf{v}_0) - (\bar{\mathbf{v}}_f - v_0)
  \end{aligned}
\end{equation}
Nonlinear dynamics of the motion can be accordingly expressed as
\begin{equation}\label{eqn:dynamic equation with nu}
  \dot{\mathbf{x}}_{\nu} =  \mathbf{f} - \mathbf{f}_0,
\end{equation}
where $\mathbf{f}$ and $\mathbf{f}_0$ are the  dynamics of the system at time $t$ and $t_0$, respectively.  Hamiltonian of the new design can then be constructed as
\begin{equation}
  \begin{aligned}
    \mathcal{H}_{\nu} & =  \boldsymbol{\lambda}_{\nu}^T\dot{\mathbf{x}}_{\nu} + \mathcal{L}_{\nu}                                                                                                                                                                                                     \\
                      & =  \boldsymbol{\lambda}_{\nu,\boldsymbol{r}} \cdot (\boldsymbol{v} - \boldsymbol{v}_0)+\boldsymbol{\lambda}_{\nu,\boldsymbol{v}} \cdot\left(\mathbf{g}(\mathbf{r})+c_1 \frac{u}{m} \boldsymbol{\alpha} - \mathbf{g}(\mathbf{r}_0)\right)-\lambda_{\nu, m} c_2 u + \nu_2 c_2 u
  \end{aligned}
\end{equation}
Corresponding costate equation of design is given as
\begin{equation}\label{eqn:costate ODE with nu}
  \dot{\boldsymbol{\lambda}}_{\nu} = -\frac{\partial \mathcal{H}_{\nu}}{\partial \mathbf{x}}=-\left(\frac{\partial \mathbf{f}}{\partial \mathbf{x}}- \frac{\partial \mathbf{f}_0}{\partial \mathbf{x}_0}\right) \boldsymbol{\lambda}_{\nu}
\end{equation}
where $\frac{\partial \mathbf{f}}{\partial \mathbf{x}}$ and $\frac{\partial \mathbf{f}_0}{\partial \mathbf{x}_0}$ are the Jacobians of the system at  time $t$ and $t_0$, respectively. Full expression of the two Jacobians are obtained using the two-body dynamics, as shown in eq.\eqref{eqn: two body costate}.

Following Pontryagin's principle, a primer-vector-based control law and switching function, similar to those in the original two-body transfer design, can be obtained as follows:
\begin{equation}
  \boldsymbol{\alpha}=-\frac{\boldsymbol{\lambda}_{\nu, \mathbf{v}}}{\lambda_{\nu, \mathbf{v}}}
\end{equation}

\begin{equation}
  \rho=1-\frac{c_1}{\nu_2 c_2 m } \lambda_{\nu, \mathbf{v}}-\frac{\lambda_{\nu, m}}{\nu_2}.   \end{equation}
Transversality condition on the mass costate, remain unchanged, and is given by $\lambda_{m,f} = 0$ since the final mass is free. Therefore, the terminal constraints for the reformulated design are summarized as follows:
\begin{equation}\label{eqn:terminal constraints with nu}
  \boldsymbol{\Psi}_\nu^f := [\delta \mathbf{r}_\nu^f, \delta \mathbf{v}_\nu^f, \lambda_{m,f}]
\end{equation}
The dynamic equations, Eq.~\eqref{eqn:dynamic equation with nu}, together with the costate ODE, Eq.~\eqref{eqn:costate ODE with nu}, and the terminal conditions in Eq.~\eqref{eqn:terminal constraints with nu}, constitute a TPBVP. Design variables of the TPBVP formulation consist of the initial costate guess, $\boldsymbol{\lambda}_{\nu,0}$, and the augmented constant, $\nu_2$.

A test is presented in the numerical simulation section to illustrate the equivalence of the two costate expressions. Initial guess solutions from the STM-based costate determination, $\boldsymbol{\lambda}_{0,\text{STM}}$, is put into the augmented TPBVP with $\nu_2$ to generate control law. The simulation results show that, for a given costate $\boldsymbol{\lambda}_{0,\text{STM}}$ and specific values of $\nu_2$, both costate TPBVP formulations can be successfully satisfied.

\section{Numerical Simulation}
\label{sec:Numerical Simulation}
A low-fidelity Earth-Mars transfer design using the proposed method is presented first. Numerical simulation results of the design are compared to those of the explicit costate-based analysis. A transfer design under high-fidelity gravity model is then presented using the proposed method. Applications of the method are further extended to a multi-revolution Earth-Venus transfer under high-fidelity models. All computations are performed in canonical units unless otherwise specified. The length unit is set to $1~\mathrm{AU}$, the time unit is defined as $t_{rcf}=l_{rcf}/v_{rcf}$, where $v_{rcf}=\sqrt{\mu/\mathrm{AU}}$ is the velocity unit, and the mass unit is set to $m_0$, the initial mass of the spacecraft.
\subsection{Indirect approach under low fidelity two-body propagator: Earth-Mars transfer}

A numerical simulation of the low-fidelity transfer design is implemented using a two-body gravity propagator. The spacecraft departs from Earth at the departure epoch 4260.62 MJD2000. An initial impulsive velocity, defined by the magnitude $v_{\infty}$, right ascension $\alpha_0$, and declination $\delta_0$, is applied at departure to initiate the transfer. The initial values of $v_{\infty}$, $\alpha_0$, and $\delta_0$ are set to 0.20 canonical units, $140^\circ$, and $3^\circ$, respectively. The initial spacecraft mass is set to $1500~\mathrm{kg}$. The thrust magnitude is $0.33~\mathrm{N}$ with a specific impulse $I_{sp}$ of $3000~\mathrm{s}$. The time of flight (ToF) is set to 474.43 days. Table~\ref{tab:Earth-Mars transfer mission parameters  and propulsion settings} lists  mission parameters and propulsion settings of the design.

\begin{table}[h!]
  \centering
  \begin{tabular}{lc}
    \toprule
    Parameter                        & Value     \\
    \midrule
    Departure epoch (MJD2000)        & $4260.62$ \\
    Initial mass (kg)                & $1500$    \\
    Time of Flight (days)            & $474.43$  \\

    $v_{\infty}$ (canonical unit)    & $0.20$    \\
    Right ascension $\alpha_0$ (deg) & $140$     \\
    Declination $\delta_0$ (deg)     & $3$       \\

    Thrust (N)                       & $0.33$    \\
    Specific impulse $I_{sp}$ (s)    & $3800$    \\
    \bottomrule
  \end{tabular}
  \caption{Earth-Mars transfer mission parameters and propulsion settings.}\label{tab:Earth-Mars transfer mission parameters  and propulsion settings}
\end{table}

The design variables for constructing the transfer consist of the initial costate guess, $\boldsymbol{\lambda}_0$, and the departure parameters, including $v_{\infty}$, $\alpha_0$, and $\delta_0$. The initial costate values are set to a uniform vector, $0.1 \times \mathbf{1}_{7}$, where $\mathbf{1}_{7}$ denotes a seven-dimensional vector of ones. Table~\ref{tab:Costate vector initialization and bounds} and Table~\ref{tab:Initial values and bounds for injection parameters} summarize the parameter settings for the initial costates and trajectory injection parameters, respectively.

\begin{table}[h!]
  \centering
  \begin{tabular}{lccc}
    \toprule
    Parameter & Initial Value                 & Lower Bound & Upper Bound \\
    \midrule
    $\boldsymbol{\lambda}_0 $
              & $ 0.1 \times \mathbf{1}_{7}$
              & $ -1.0 \times \mathbf{1}_{7}$
              & $ 1.0 \times \mathbf{1}_{7}$                              \\
    \bottomrule
  \end{tabular}
  \caption{Costate vector initialization and bounds} \label{tab:Costate vector initialization and bounds}
\end{table}

\begin{table}[h!]
  \centering
  \begin{tabular}{lccc}
    \toprule
    Parameter & Initial Value & Lower Bound & Upper Bound \\
    \midrule
    $ v_{\infty}$
              & $0.20$
              & $0.20 - 0.05$
              & $0.20 + 0.05$                             \\

    $ \alpha_0 $ (deg)
              & $140$
              & $140 - 18$
              & $140 + 18$                                \\

    $ \delta_0$  (deg)
              & $3$
              & $3 - 5$
              & $3 + 5$                                   \\
    \bottomrule
  \end{tabular}
  \caption{Initial values and bounds for injection parameters}\label{tab:Initial values and bounds for injection parameters}
\end{table}

The terminal boundary conditions, given by the target ephemeris of Mars at the arrival epoch, are computed using the ephemeris propagation function and subsequently normalized in canonical units. The transversality condition for the mass costate, $\lambda_{m,f}$, is set to zero at the final time. Table~\ref{tab:Terminal conditions for the Earth--Mars transfer in canonical units} lists the terminal conditions for the TPBVP of the transfer design.

\begin{table}[h!]
  \centering
  \setlength{\tabcolsep}{0pt}
  \begin{tabular*}{\textwidth}{@{\extracolsep{\fill}}lccc}
    \toprule
    Parameter & $\bar{\mathbf{r}}_f$ & $\bar{\mathbf{v}}_f$ & $\bar{\lambda}_{m,f}$ \\
    \midrule
    Value &
    $[0.94763, -1.0174, -0.044592]$ &
    $[0.62590, 0.62425, -0.0023554]$ &
    $0$ \\
    \bottomrule
  \end{tabular*}
  \caption{Terminal conditions for the Earth--Mars transfer} \label{tab:Terminal conditions for the Earth--Mars transfer in canonical units}
\end{table}

The system equations of motion and the STM-based costate equations are propagated using an explicit Runge--Kutta solver, \texttt{ode45}. The numerical integration tolerances, \texttt{AbsTol} and \texttt{RelTol}, are set to $10^{-10}$ and $10^{-11}$, respectively. A standard Sequential Quadratic Programming (SQP) optimizer is used to search for the optimal initial costates. The maximum number of function evaluations is limited to $2000$.

The solution converges with a total residual error on the order of $10^{-6}$. The terminal constraint errors, including the position and velocity mismatches and the transversality condition on the mass costate, are listed in Table~\ref{tab:Terminal constraint errors for the rendezvous solution}. The optimized decision vector is given in Table~\ref{tab:Optimized initial costates for the low-fidelity Earth--Mars rendezvous solution} and Table~\ref{tab:Optimized departure parameters}. The solution yields a final mass of $1341.9~\mathrm{kg}$. The resulting mission performance metrics are summarized in Table~\ref{tab:Resulting mission performance metrics}.

\begin{table}[h!]
  \centering
  \begin{tabular*}{\textwidth}{@{\extracolsep{\fill}}lcccc}
    \toprule
    Quantity & $\|\Delta \mathbf{r}_f\|$ & $\|\Delta \mathbf{v}_f\|$ & $\lambda_{m,f}$ & Total error norm \\
    \midrule
    Error norm & $2.89e-6$ & $3.31e-6$ & $5.3608e-8$ & $4.3918e-6$ \\
    \bottomrule
  \end{tabular*}
  \caption{Terminal Constraint Error Norms for the Two-Body Gravity Model}
  \label{tab:Terminal constraint errors for the rendezvous solution}
\end{table}

\begin{table}[h!]
  \centering
  \renewcommand{\arraystretch}{1.15}
  \begin{tabular}{lc}
    \toprule
    Costate           & Value       \\
    \midrule
    $\lambda_{x,0}$   & $-0.14019$  \\
    $\lambda_{y,0}$   & $-0.75479$  \\
    $\lambda_{z,0}$   & $-0.01694$  \\
    $\lambda_{v_x,0}$ & $-0.48182$  \\
    $\lambda_{v_y,0}$ & $0.47409$   \\
    $\lambda_{v_z,0}$ & $0.17477$   \\
    $\lambda_{m,0}$   & $-0.054846$ \\
    \bottomrule
  \end{tabular}
  \caption{Optimized Initial Costates for the Two-Body Gravity Model}\label{tab:Optimized initial costates for the low-fidelity Earth--Mars rendezvous solution}
\end{table}

\begin{table}[h!]
  \centering
  \renewcommand{\arraystretch}{1.15}
  \begin{tabular}{lccc}
    \toprule
    Parameter & $ v_{\infty}$ & $\alpha_0$ (rad) & $\delta_0$ (rad) \\
    \midrule
    Value     & $0.238467$    & $2.49321$        & $0.0550926$      \\
    \bottomrule
  \end{tabular}
  \caption{Optimized Departure Parameters for the Two-Body Gravity Model}\label{tab:Optimized departure parameters}
\end{table}

\begin{table}[h!]
  \centering
  \renewcommand{\arraystretch}{1.15}
  \begin{tabular}{lccc}
    \toprule
    Quantity & Time of Flight (days) & Final Mass & Total Error Norm \\
    \midrule
    Value    & $474.43$              & $1341.9$   & $4.3918e-6$      \\
    \bottomrule
  \end{tabular}
  \caption{Resulting Mission Performance Metrics for the Two-Body Gravity Model}\label{tab:Resulting mission performance metrics}
\end{table}

\begin{figure}[h]
  \centering
  \includegraphics[width=0.7\linewidth]{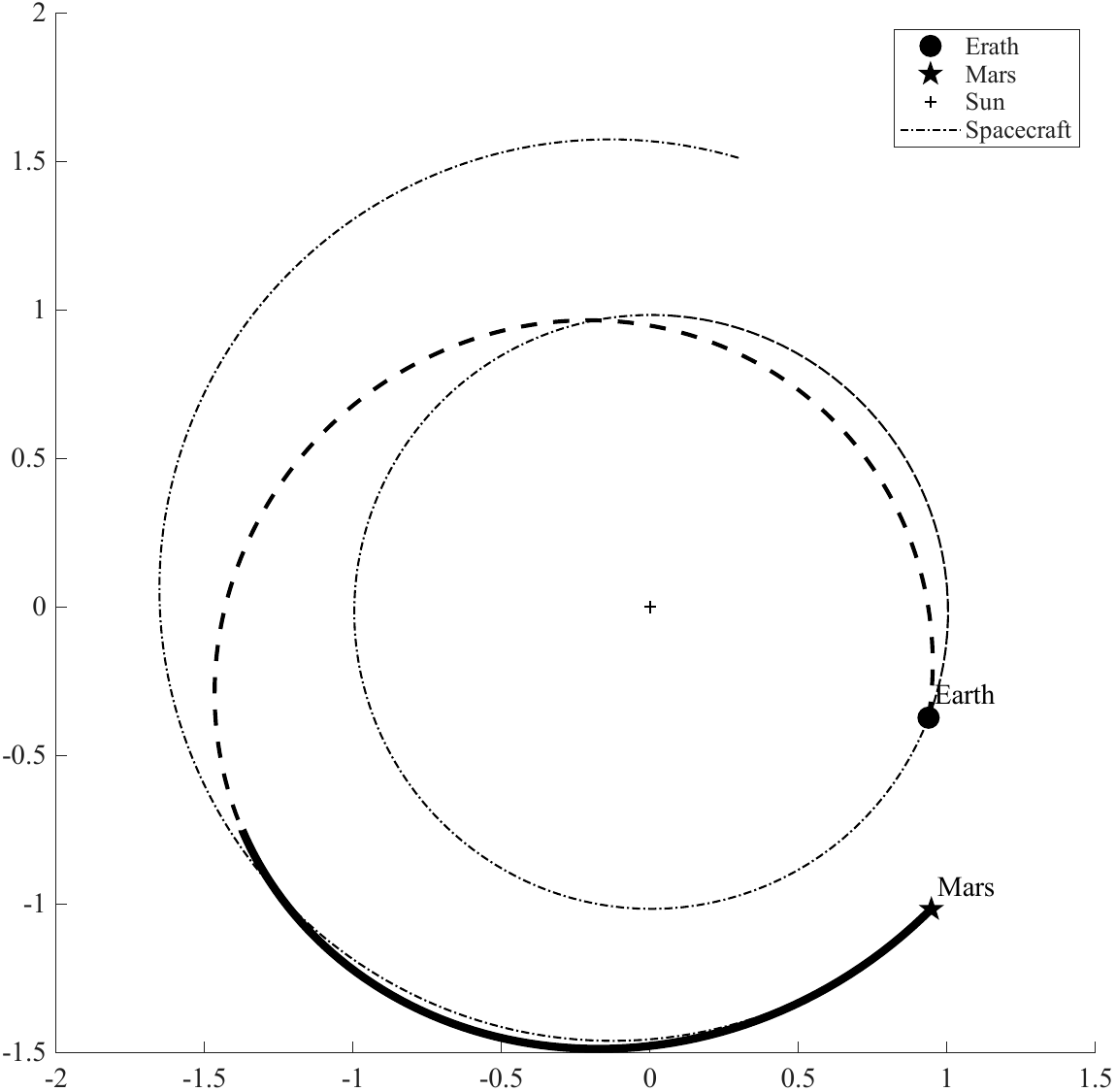}
  \caption{Earth-Mars Trajectory Design in the Two-Body Gravity Model}
  \label{fig:low_fidelity_earth_mars_trajectory}
\end{figure}

\begin{figure}[h]
  \centering
  \includegraphics[width=0.7\linewidth]{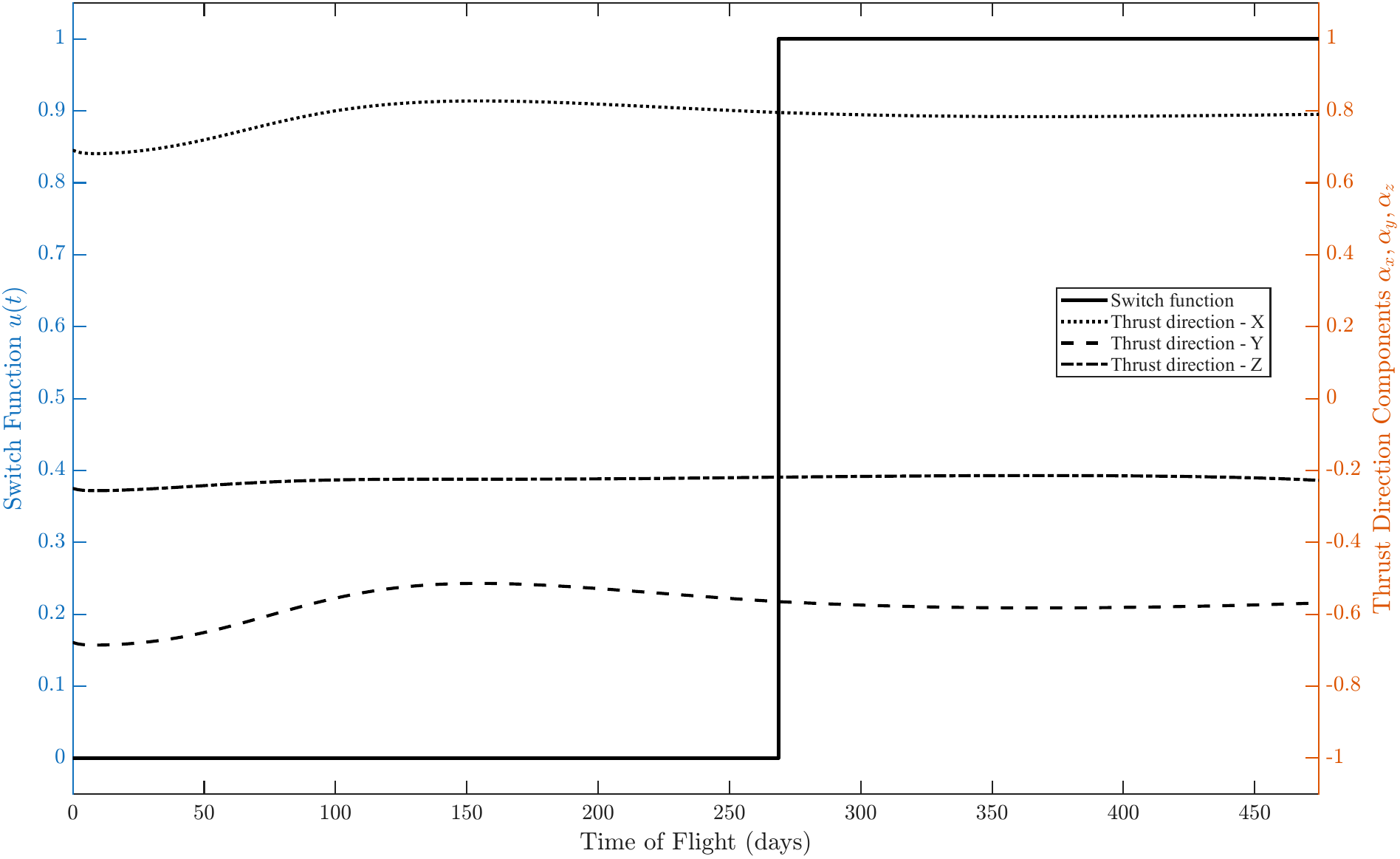}
  \caption{Thrust Profile of the Earth-Mars Transfer in the Two-Body Gravity Model}
  \label{fig:low_fidelity_earth_mars_thrust}
\end{figure}

The trajectory and thrust profile of the transfer are shown in Figure~\ref{fig:low_fidelity_earth_mars_trajectory} and Figure~\ref{fig:low_fidelity_earth_mars_thrust}, respectively. The active thrusting and coasting phases of the trajectory are distinguished using different line styles.

\subsection{Equivalence Verification to the Explicit Costate Equation in the Two-Body Gravity Model}

A numerical test is conducted to demonstrate the equivalence between the proposed initial costate determination and the explicit costate-equation approach. Design parameters of the explicit costate-based formulation are set to the augmented objective parameter $\nu_2$ and the initial costates $\boldsymbol{\lambda}_0$. The optimal initial costates obtained from the STM-based approach, $\boldsymbol{\lambda}_{0,\mathrm{STM}}$, are used as the starting point for $\boldsymbol{\lambda}_0$. The departure parameters for the test, listed in Table~\ref{tab:Departure parameters for the costate-equation equivalence test}, are set to the same values as those of the STM-based solution in Table~\ref{tab:Optimized departure parameters}. The trajectory ODE integrator and error tolerance settings are also kept identical to those used in the STM-based design. Table~\ref{tab:initial_costate_augmented_parameter_test} and Table~\ref{tab:Departure parameters for the costate-equation equivalence test} summarize the parameter settings of the verification test.

\begin{table}[h!]
  \centering
  \renewcommand{\arraystretch}{1.15}
  \begin{tabular}{lccc}
    \toprule
    Parameter                  & Initial value                           & Lower bound                                  & Upper bound                                  \\
    \midrule
    $\boldsymbol{\lambda}_{0}$ & $\boldsymbol{\lambda}_{0,\mathrm{STM}}$ & $\boldsymbol{\lambda}_{0,\mathrm{STM}}-0.08$ & $\boldsymbol{\lambda}_{0,\mathrm{STM}}+0.08$ \\
    $\nu_2$                    & $0.95$                                  & $0.10$                                       & $1.50$                                       \\
    \bottomrule
  \end{tabular}
  \caption{Initial Parameter Settings for the Costate-Equation Equivalence Test}
  \label{tab:initial_costate_augmented_parameter_test}
\end{table}

\begin{table}[h!]
  \centering
  \renewcommand{\arraystretch}{1.15}
  \begin{tabular}{lccc}
    \toprule
    Parameter & $ v_{\infty}$ & $\alpha_0$ (rad) & $\delta_0$ (rad) \\
    \midrule
    Value     & $0.238467$    & $2.49321$        & $0.0550926$      \\
    \bottomrule
  \end{tabular}
  \caption{Departure Parameters for the Costate-Equation Equivalence Test}\label{tab:Departure parameters for the costate-equation equivalence test}
\end{table}

The explicit costate equation in Eq.~\eqref{eqn:costate with Jacobian}, together with the two-body dynamics, is implemented. The gradient-based SQP optimizer is used to search for the solution. Table~\ref{tab:Terminal constraint error norms: explicit costate equation} lists the terminal constraint error norms of the resulting trajectory. Table~\ref{tab:lambda_comparison_stm_costate_low_fidelity} presents the numerical solution and compares it with the STM-based solution. The optimized augmented parameter $\nu_2$ is also listed. The explicit costate-equation solution remains close to the STM-based solution for most components. Table~\ref{tab:stm_explicit_costate_performance_comparison} compares the final mass and residual norm of the two solutions. The two approaches produce nearly identical final masses with similar final residual norms.

\begin{table}[h!]
  \centering
  \begin{tabular*}{\textwidth}{@{\extracolsep{\fill}}lcccc}
    \toprule
    Quantity & $\|\Delta \mathbf{r}_f\|$ & $\|\Delta \mathbf{v}_f\|$ & $\lambda_{m,f}$ & Total error norm \\
    \midrule
    Error norm &  1.0956e-04 & 4.0575e-04 & 6.3179e-08 & 4.2028e-04 \\
    \bottomrule
  \end{tabular*}
  \caption{Terminal Constraint Error Norms for the Explicit Costate Equation}
  \label{tab:Terminal constraint error norms: explicit costate equation}
\end{table}

\begin{table}[h!]
  \centering
  \renewcommand{\arraystretch}{1.15}
  \begin{tabular}{lccc}
    \toprule
    Solution          & STM       & Explicit & Difference \\
    \midrule
    $\lambda_{x,0}$   & -0.14019  & -0.14019 & 0          \\
    $\lambda_{y,0}$   & -0.75479  & -0.75479 & 0          \\
    $\lambda_{z,0}$   & -0.01694  & -0.01694 & 0          \\
    $\lambda_{v_x,0}$ & -0.48182  & -0.54576 & -0.06394   \\
    $\lambda_{v_y,0}$ & 0.47409   & 0.40452  & -0.06957   \\
    $\lambda_{v_z,0}$ & 0.17477   & 0.14932  & -0.02545   \\
    $\lambda_{m,0}$   & -0.054846 & -0.0863  & -0.031454  \\
    $\nu_2$           & $-$       & 0.88093  & $-$        \\
    \bottomrule
  \end{tabular}
  \caption{Costate Solution for the Explicit Costate Equation}
  \label{tab:lambda_comparison_stm_costate_low_fidelity}
\end{table}

\begin{table}[h!]
  \centering
  \renewcommand{\arraystretch}{1.15}
  \small
  \begin{tabular*}{\textwidth}{@{\extracolsep{\fill}}lccc}
    \toprule
    Method & ToF (days) & Final mass (kg) & Final residual norm \\
    \midrule
    STM method & 474.43 & 1341.9 & 4.3918e-6 \\
    Explicit costate & 474.43 & 1342.7 & 4.2028e-4 \\
    \bottomrule
  \end{tabular*}
  \caption{Performance Comparison Between STM-Based and Explicit Costate Solutions}
  \label{tab:stm_explicit_costate_performance_comparison}
\end{table}

\subsection{Earth-Mars Transfer Design Under a High-Fidelity Model}

An Earth-Mars transfer design under high fidelity dynamic model is conducted. Impacts of central Sun gravity, Jupiter perturbation (dominant third body), solar radiation pressure, solar J2 and first-order post-Newtonian relativistic correction, etc., are considered (Table \ref{tab:high_fidelity_model_options}). The total orbital acceleration and propagator computations are configured according to Eq.~\eqref{eqn:high_fidelity_acc} and Eq.~\eqref{eqn:high_fidelity_propagator}. The design scenario parameters, including the departure epoch, initial mass, and propulsion settings, are set to the same values as those used in the low-fidelity gravity-propagator model (Table~\ref{tab:Earth-Mars transfer mission parameters  and propulsion settings}). Configuration and settings of the initial costates, departure parameters, and terminal constraints are also from the low-fidelity design. Table~\ref{tab:low_high_config_settings_compare} summarizes and compares the nonlinear optimizer configuration and force-model settings for the two transfer designs.

\begin{table}[h!]
  \centering
  \renewcommand{\arraystretch}{1.15}
  \small
  \begin{tabular*}{\textwidth}{@{\extracolsep{\fill}} l l l}
    \toprule
    Model option & Included effect & Parameter / reference \\
    \midrule
    Central Sun gravity & Point-mass solar attraction & eq.\eqref{eqn:high_fidelity_acc} \\
    Jupiter perturbation & Dominant third-body acceleration & Jupiter ephemeris, $\mu_J$ \\
    Solar radiation pressure & SRP acceleration & $\beta = 1\times10^{-7}$ \\
    Solar $J_2$ & Sun quadrupole gravity perturbation & $J_2 = 2.2\times10^{-7}$ \\
    Relativistic correction & First-order post-Newtonian term & $c = 173.1446327$ \\
    \bottomrule
  \end{tabular*}
  \caption{High-Fidelity Model Options}
  \label{tab:high_fidelity_model_options}
\end{table}

The equations of motion and the STM-based costate equations in Eq.~\eqref{eqn:STM-based costate equation} are integrated using a variable-step, variable-order Adams solver, \texttt{ode113}. Absolute and relative tolerances of the numerical integrator, \texttt{AbsTol} and \texttt{RelTol}, are set to $10^{-12}$ and $10^{-11}$, respectively. The trajectory optimization problem is solved using a gradient-based SQP solver.
\begin{table}[h!]
  \centering
  \begin{tabular}{ccc}
    \toprule
    Setting / Feature & Low-Fidelity                     & High-Fidelity                                    \\
    \midrule
    Initial costates  & $0.1\times\mathbf{1}_7$          & $0.1\times\mathbf{1}_7$                          \\
    Costate bounds    & $[-\mathbf{1}_7,\,\mathbf{1}_7]$ & $[-\mathbf{1}_7,\,\mathbf{1}_7]$                 \\
    $ v_{\infty}$     & $[0.2 \pm 0.05]$                 & $[0.2 \pm 0.05]$                                 \\
    $\alpha_0 $       & $[140.0\pm18]$                   & $[140.0\pm18]$                                   \\
    $\delta_0 $       & $[3.0\pm5.0]$                    & $[3.0\pm5.0]$                                    \\
    ODE Solver        & $\texttt{ode45}$                 & $\texttt{ode113}$                                \\
    Force Model       & Two-body gravity                 & High-fidelity (eq.\eqref{eqn:high_fidelity_acc}) \\
    \bottomrule
  \end{tabular}
  \caption{Numerical Configuration Settings for the Earth-Mars Transfer Designs}
  \label{tab:low_high_config_settings_compare}
\end{table}

The optimization converges after approximately 2000 function evaluations. Table~\ref{tab:departure parameters of the high-fidelity solution} summarizes the optimized departure parameters under the low- and high-fidelity models, and Table~\ref{tab:costate_compare} compares the corresponding initial costate solutions. The high-fidelity model leads to slight changes in the optimal departure velocity magnitude and direction. Several costate components, including $\lambda_{0,x}$, $\lambda_{0,z}$, and $\lambda_{0,vx}$, are also altered, reflecting effects of the additional perturbations included in the heliocentric dynamics.
\begin{table}[h!]
  \centering
  \renewcommand{\arraystretch}{1.15}
  \begin{tabular}{lcc}
    \toprule
    Departure parameter & Low Fidelity & High Fidelity \\
    \midrule
    $ v_{\infty}$       & $0.23846$    & $0.229477$    \\
    $\alpha_0$ (rad)    & $2.49321$    & $2.474037$    \\
    $\delta_0$ (rad)    & $0.05509$    & $0.051139$    \\
    \bottomrule
  \end{tabular}
  \caption{Optimized Departure Parameters for the Low- and High-Fidelity Gravity Models}
  \label{tab:departure parameters of the high-fidelity solution}
\end{table}

\begin{table}[h!]
  \centering
  \begin{tabular}{lcc}
    \toprule
    Costate          & Low Fidelity & High Fidelity \\
    \midrule
    $\lambda_{0,x}$  & $-0.14019$   & $-0.28526$    \\
    $\lambda_{0,y}$  & $-0.75479$   & $-0.66883$    \\
    $\lambda_{0,z}$  & $-0.01694$   & $0.18613$     \\
    $\lambda_{0,vx}$ & $-0.48182$   & $-0.55339$    \\
    $\lambda_{0,vy}$ & $0.47409$    & $0.45125$     \\
    $\lambda_{0,vz}$ & $0.17477$    & $0.17403$     \\
    $\lambda_{0,m}$  & $-0.054846$  & $0.010007$    \\
    \bottomrule
  \end{tabular}
  \caption{Costate Solutions for the Earth-Mars Rendezvous}
  \label{tab:costate_compare}
\end{table}

\begin{figure}[!h]
  \centering
  \includegraphics[width=0.7\linewidth]{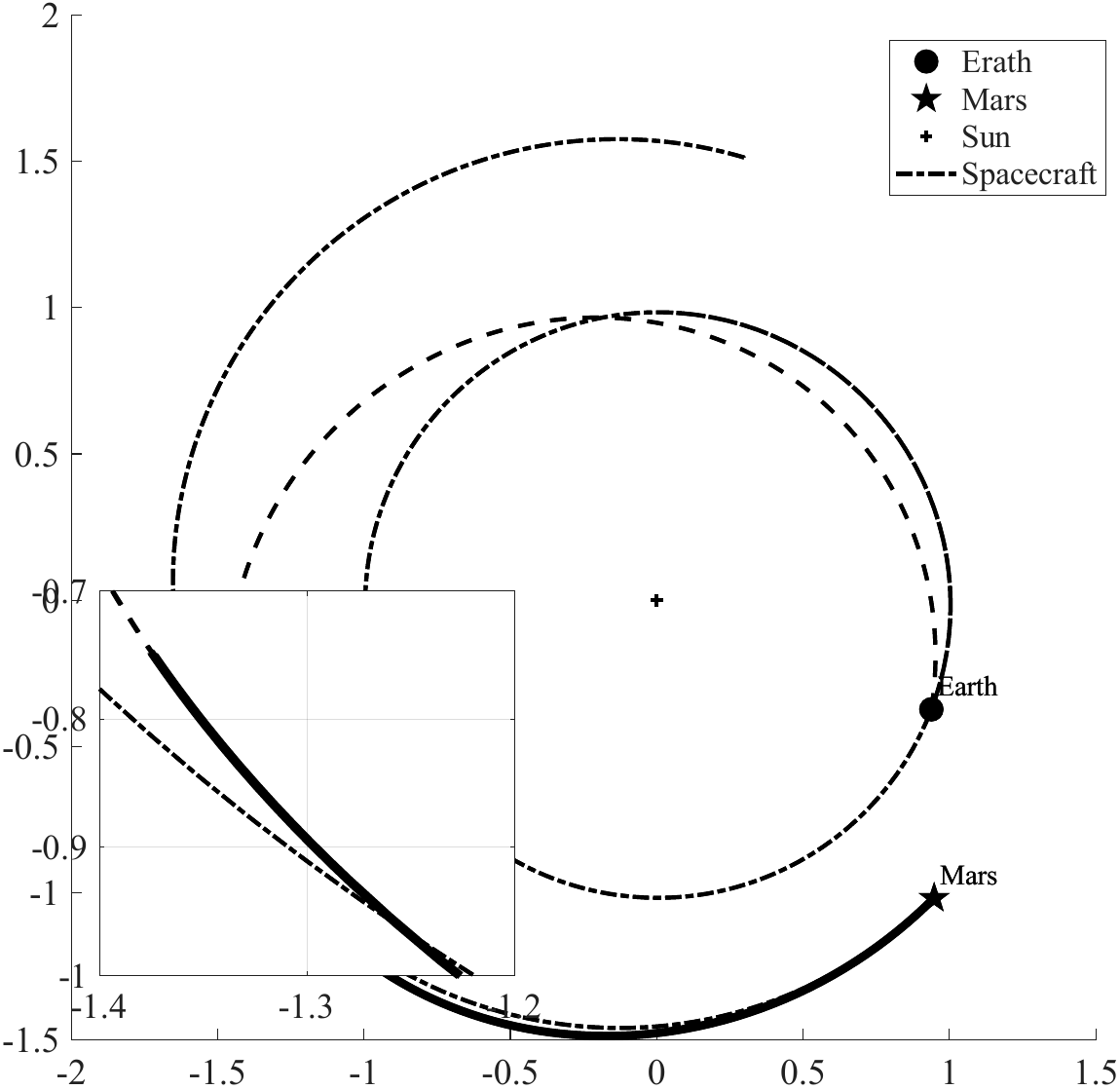}
  \caption{Comparison of Low- and High-Fidelity Earth-Mars Rendezvous Trajectory Designs}
  \label{fig:low_high_fidelity_trajectory_comparison}
\end{figure}

\begin{figure}[!h]
  \centering
  \includegraphics[width=0.7\linewidth]{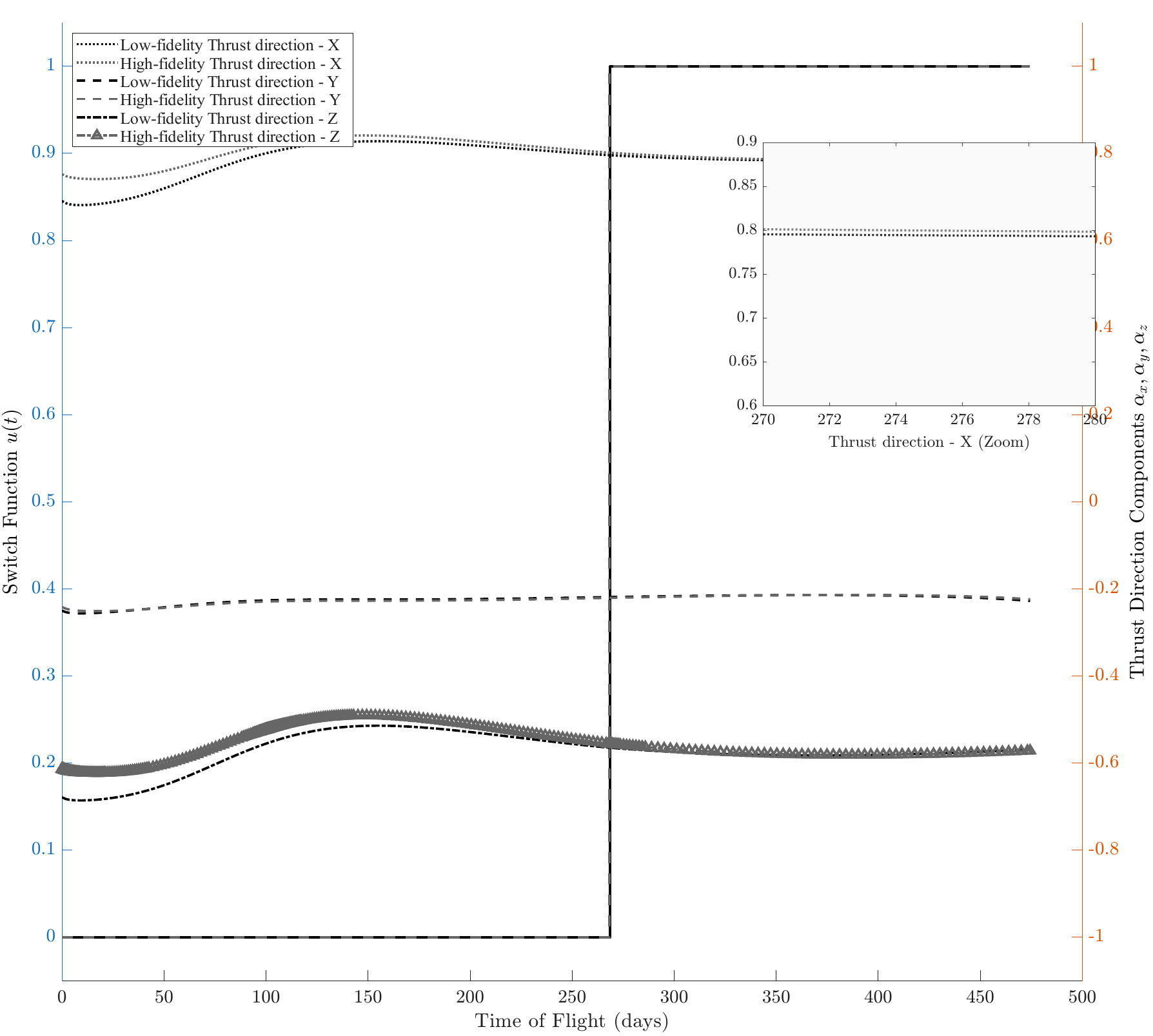}
  \caption{Comparison of Low- and High-Fidelity Earth-Mars Rendezvous Thrust Profiles}
  \label{fig:low_high_fidelity_thrust_comparison}
\end{figure}

Figures~\ref{fig:low_high_fidelity_trajectory_comparison} and \ref{fig:low_high_fidelity_thrust_comparison} present and compare the trajectories and thrust profiles of the two designs. The trajectories are nearly identical. Both contain a single thrust arc starting at nearly the same point. The thrust profile of the high-fidelity solution shows slightly different thrust directions from its low-fidelity counterpart to account for the included higher-order perturbations.

Table~\ref{tab:terminal_compare} lists terminal constraint errors of the two designs. Both models achieve similar residual levels, with total error norms on the order of $10^{-5}$ and $10^{-6}$ in canonical units, respectively, and both successfully satisfy the rendezvous conditions. Table~\ref{tab:flight_mass_residual_comparison} presents the performance metrics of the low- and high-fidelity rendezvous solutions. The final spacecraft mass for the high-fidelity design is $1350.9~\mathrm{kg}$. The value is slightly higher than that of the low-fidelity counterpart, while the same time of flight is preserved.

\begin{table}[h!]
  \centering
  \begin{tabular}{lcc}
    \toprule
    Constraint Norm           & Low Fidelity & High Fidelity \\
    \midrule
    $\|\Delta \mathbf{r}_f\|$ & 2.89e-6      & 1.06e-5       \\
    $\|\Delta \mathbf{v}_f\|$ & 3.31e-6      & 1.72e-5       \\
    $|\lambda_m(t_f)|$        & 5.36e-8      & 1.80e-5       \\
    Total Error Norm          & 4.39e-6      & 2.72e-5       \\
    \bottomrule
  \end{tabular}
  \caption{Comparison of Terminal Constraint Norms Between the Low- and High-Fidelity Solutions}
  \label{tab:terminal_compare}
\end{table}

\begin{table}[h!]
  \centering
  \renewcommand{\arraystretch}{1.15}
  \begin{tabular}{lcc}
    \toprule
    Quantity              & Low Fidelity & High Fidelity \\
    \midrule
    Time of Flight (days) & 474.43       & 474.43        \\
    Final Mass (kg)       & 1341.9       & 1350.9        \\
    Final Residual Norm   & 4.3918e-6    & 2.7213e-5     \\
    \bottomrule
  \end{tabular}
  \caption{Performance Metrics for the Low- and High-Fidelity Earth-Mars Rendezvous Solutions}
  \label{tab:flight_mass_residual_comparison}
\end{table}

\subsection{Application to High Fidelity Earth-Venus Transfer Design}

A numerical simulation of the multiple revolution Earth-Venus transfer design is implemented. The transfer consists of multiple coast-thrust arcs, with a total time of flight of 1000 days. The same high-fidelity force model numerical settings to those in the Earth-Mars transfer design are used.

Propulsion settings and initial mass of the transfer are set to the same values as those used in the Earth-Mars transfer design, too. The thrust magnitude, $I_{sp}$ and initial mass are set to 0.33N, 3800s and 1500kg, respectively. Table \ref{tab:EV_boundary_conditions} shows departure and terminal conditions of the design.   Initial value of design variables, $\boldsymbol{\lambda}_0 $, is set to $ 0.3 \times \mathbf{1}_{7}$. Table \ref{tab:EV_lambda_bounds}  summarizes parameter settings of design varaiables of the initial costate.

\begin{table}[t]
  \centering
  \begin{tabular}{lcc   }
    \toprule
    {State} & {Departure}
             & {Terminal}                \\
    \midrule
    $x$      & 0.9708322    & -0.3277178 \\
    $y$      & 0.2375844    & 0.6389172  \\
    $z$      & -1.671055e-6 & 0.02765929 \\
    $v_x$    & -0.2543859   & -1.050138  \\
    $v_y$    & 0.9679734    & -0.5431309 \\
    $v_z$    & 1.502959e-5  & 0.05317731 \\
    \bottomrule
  \end{tabular}
  \caption{Departure and Terminal Conditions for the Earth-Venus Transfer}
  \label{tab:EV_boundary_conditions}
\end{table}

\begin{table}[h!]
  \centering
  \begin{tabular}{lccc}
    \toprule
    Parameter                    & Initial Value            & Lower Bound               & Upper Bound              \\
    \midrule
    $\boldsymbol{\lambda}_{0,r}$ & $0.3\times \mathbf{1}_3$ & $-1.0\times \mathbf{1}_3$ & $1.0\times \mathbf{1}_3$ \\
    $\boldsymbol{\lambda}_{0,v}$ & $0.3\times \mathbf{1}_3$ & $-1.0\times \mathbf{1}_3$ & $1.0\times \mathbf{1}_3$ \\
    $\lambda_{0,m}$              & $0.3$                    & $0.0$                     & $1.0$                    \\
    \bottomrule
  \end{tabular}
  \caption{Initial Values and Bounds for the Costate Variables}
  \label{tab:EV_lambda_bounds}
\end{table}

\begin{figure}[h!]
  \centering
  \includegraphics[width=0.7\linewidth]{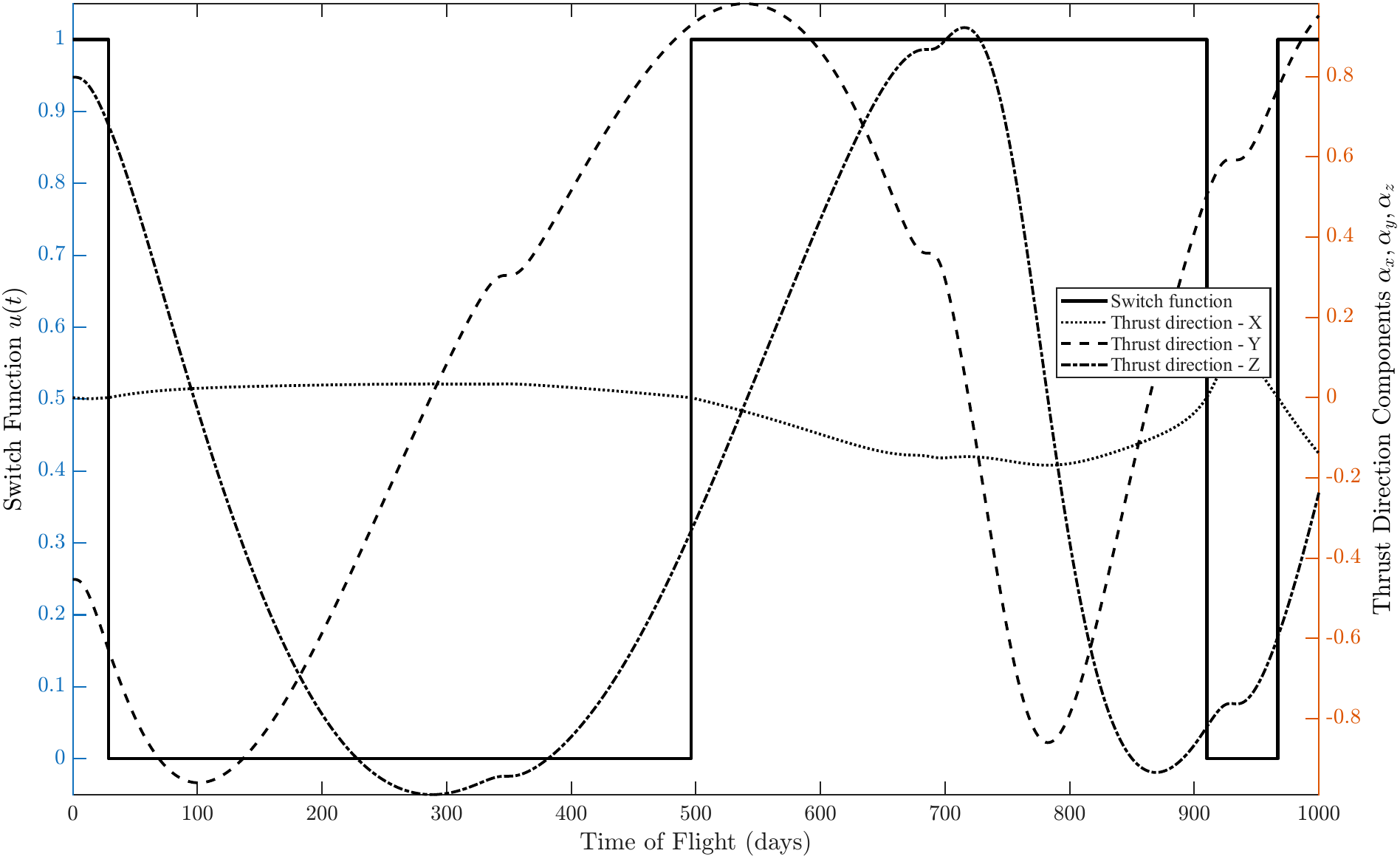}
  \caption{Earth-Venus thrust profile}
  \label{fig:earth_VENUS_thrust}
\end{figure}

\begin{figure}[h!]
  \centering
  \includegraphics[width=0.7\linewidth]{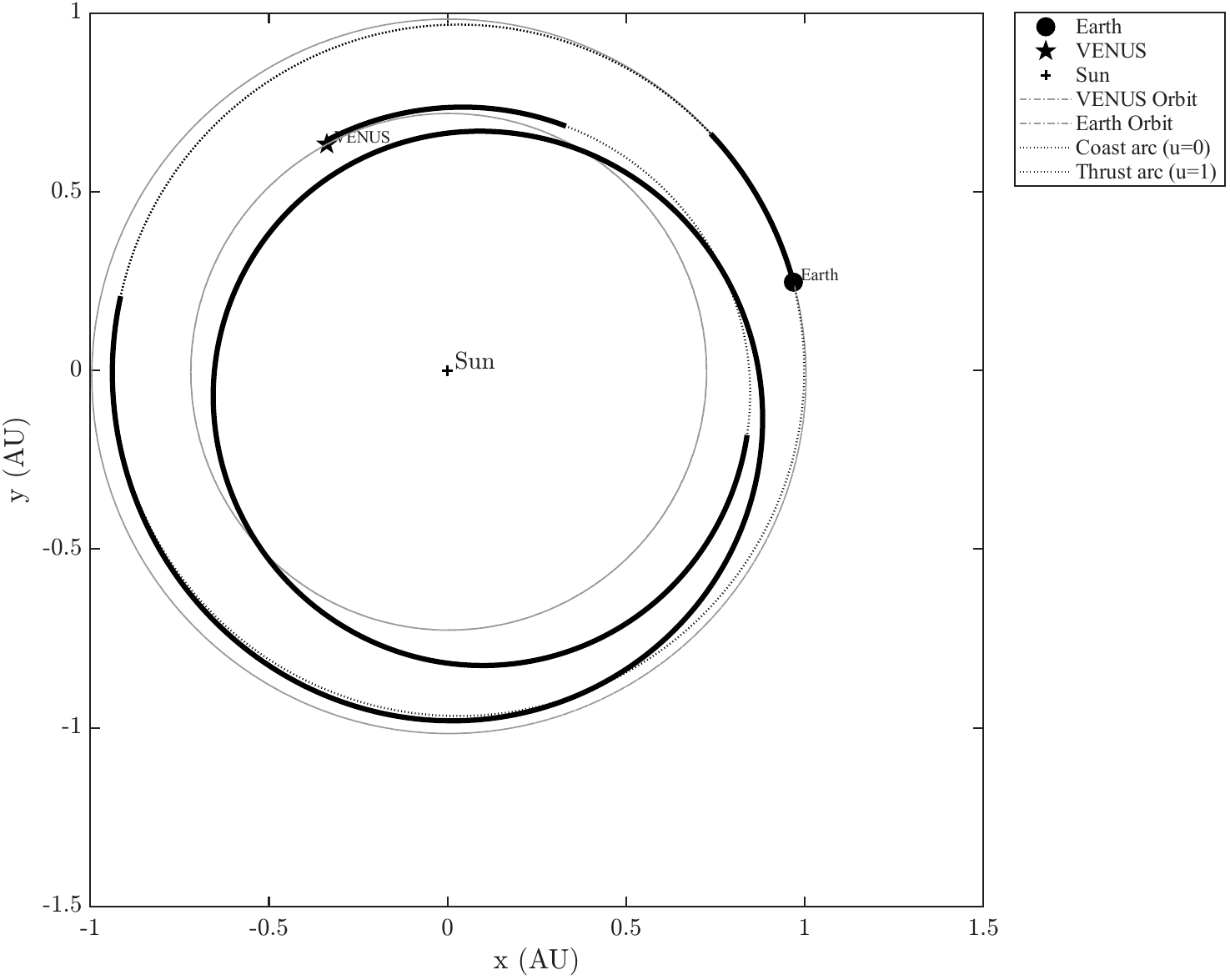}
  \caption{Earth-Venus trajectory design: High-fidelity model}
  \label{fig:earth_VENUS_trajectory}
\end{figure}

The variable‑order, multistep \texttt{ode113} integrator is implemented for numerical integration of the trajectory and costate equations.  Absolute and relative error tolerances of the ODE solver are set to 1e-11 and 1e-12, respectively. The optimizer option is set to the interior-point algorithm, with the maximum number of function evaluations limited to $6000$. Figures \ref{fig:earth_VENUS_thrust} and \ref{fig:earth_VENUS_trajectory} present thrust profile and trajectory of the  transfer design. The multiple coast-thrust arcs vary according to the switching function. The optimal initial costates and departure parameters are listed in Table~\ref{tab:EV_optimal_lambda0}. Table~\ref{tab:EV_terminal_shooting} shows terminal constraint errors. Performance metrics of the design are summarized in Table~\ref{tab:EV_performance_metrics}. The solution yields a final mass of $1136.25~\mathrm{kg}$ with a final error norm of $0.07731$.

The design converges to a rendezvous solution with similar performance metrics and terminal constraint errors as those obtained in the Earth-Mars cases.  Applicability and robustness of the proposed indirect optimization method are demonstrated through  different mission scenarios and target bodies.

\begin{table}[h!]
  \centering
  \begin{tabular}{lc}
    \toprule
    Component         & Value      \\
    \midrule
    $\lambda_{0,x}$   & 0.42905083 \\
    $\lambda_{0,y}$   & 0.28222769 \\
    $\lambda_{0,z}$   & 0.53702772 \\
    $\lambda_{0,v_x}$ & -0.2535146 \\
    $\lambda_{0,v_y}$ & 0.46718312 \\
    $\lambda_{0,v_z}$ & 0.35312778 \\
    $\lambda_{0,m}$   & 0.20158164 \\
    \bottomrule
  \end{tabular}
  \caption{Optimal Initial Costate Vector for the Earth-Venus Transfer}
  \label{tab:EV_optimal_lambda0}
\end{table}

\begin{table}[h!]
  \centering
  \begin{tabular}{cc}
    \toprule
    Parameter        & Value        \\
    \midrule
    $\Delta x_{f}$   & -6.691599e-3 \\
    $\Delta y_{f}$   & 7.378535e-3  \\
    $\Delta z_{f}$   & -2.988433e-2 \\
    $\Delta v_{f,x}$ & -5.851071e-2 \\
    $\Delta v_{f,y}$ & -2.408828e-2 \\
    $\Delta v_{f,z}$ & -2.426931e-2 \\
    $\lambda_{f,m}$  & 1.980621e-2  \\
    \bottomrule
  \end{tabular}
  \caption{Terminal Shooting Constraint Errors of the Earth-Venus Transfer}
  \label{tab:EV_terminal_shooting}
\end{table}

\begin{table}[h!]
  \centering
  \renewcommand{\arraystretch}{1.15}
  \small
  \begin{tabular*}{\textwidth}{@{\extracolsep{\fill}}lccc}
    \toprule
    Quantity & Time of Flight (days) & Final Mass & Total Error Norm \\
    \midrule
    Value & 1000       & 1136.25         & 0.07731          \\
    \bottomrule
  \end{tabular*}
  \caption{Performance Metrics for the Earth-Venus Transfer Solution}
  \label{tab:EV_performance_metrics}
\end{table}

\section{Conclusion}\label{sec:conclusion}

Conventional indirect approach to the space trajectory design requires explicit costate expressions, of which only comparatively simple gravity models, e.g., two or three-body dynamics are available. With the new developed method, high-order nonlinear impacts, such as the solar pressure, Jupiter three-body impacts, Sun J2 oblateness, and relativity effects, etc., can be  incorporated and treated as a black-box in the design optimization. A standard gradient or interior-point based optimizer can be used to directly construct the optimal high fidelity transfer. The proposed design attains high‑order rendezvous accuracy when evaluated within the high‑fidelity dynamical environment. Solution of the method is verified through comparison with the explicit costate equation solution.

To accurately capture switching function in the design optimization, a stringent error tolerance, such as 1.0e-10 for both relative and total error tolerances, should be set with a high order numerical integrator.  Such numerical settings may introduce additional computational cost for medium- and large-scale mission designs involving multiple spacecraft and long-duration trajectories.  Further improvement of the proposed optimization approach is underway. Preliminary results from a variant of the method appear promising, showing greater robustness and lower requirements on computational settings. The method is also being extended to more challenging designs with chaotic characteristics,  such as low-thrust tether-based orbital debris removal. Space designs for large-scale and sequential low-thrust space operations will also be considered in future work.

\iffastcompile
  \IfFileExists{main-elsarticle-template-num-0423-3.bbl}{\input{main-elsarticle-template-num-0423-3.bbl}}{}
\else
  \bibliographystyle{elsarticle-num}
  \bibliography{reference}
\fi



\end{document}